\documentclass{article}

\usepackage{latexsym}
\usepackage{amsfonts,amssymb}
\usepackage{amsmath}

\newtheorem{theorem}{Theorem}[section]
\newtheorem{lemma}[theorem]{Lemma}
\newtheorem{corollary}[theorem]{Corollary}
\newtheorem{remark}[theorem]{Remark}

\newtheorem{proposition}[theorem]{Proposition}

\def\Mod{\mbox{\rm{Mod}}}

\def\Diff{\mbox{\rm{Diff}}}
\def\Stab{\mbox{\rm{Stab}}}
\def\SL{\mbox{\rm{SL}}}

\def\Aut{\mbox{\rm{Aut}}}

\def\Hom{\mbox{\rm{Hom}}}

\def\Triv{\mbox{\rm{Triv}}}
\def\Out{\mbox{\rm{Out}}}
\def\Ga{\Gamma}
\def\C{\bf C}
\def\La{\Lambda}
\def\Z{\bf Z}
\def\ga{\gamma}
\def\Sym{\mbox{\rm{Sym}}}
\def\rank{\mbox{\rm{rank}}}

\begin{document}


\title{Separable subgroups of mapping class groups}
\author{Christopher J. Leininger \thanks{Supported by an N.S.F. postdoctoral fellowship.}\\ D. B. McReynolds\thanks{Supported by a Continuing Education fellowship.}}
\maketitle


\begin{abstract}
We investigate the separability of several well known classes of 
subgroups of the mapping class group of a surface.
\end{abstract}


\bibliographystyle{amsplain}

\section{Introduction and main results}\label{section1}

\subsection{Introduction}

A subgroup $H < G$ is said to be \emph{separable} in $G$ if it can be expressed as an intersection of finite index subgroups of $G$. If the trivial subgroup $\{1 \}$ is separable, $G$ is \emph{residually finite}.
More generally, if every finitely generated subgroup of $G$ is separable, $G$ is \emph{subgroup separable} (or \emph{LERF}).

Subgroup separability has been an important tool in geometry; for example it often permits the lifting of a $\pi_1$--injective immersion to an embedding in a finite cover \cite{Scott78} (see in addition \cite{Bergeron00}, \cite{LongReid02,LongReid04B}, \cite{Long87,Long88}, and \cite{McReynolds04A}). Algebraically, it can also be viewed as an indication of an abundance of finite index subgroups and a rich interaction of these subgroups with the finitely generated ones. This powerful property is generally difficult to establish and the class of groups known to be subgroup separable is small. It is a theorem of M. Hall \cite{Hall49} that free groups are subgroup separable. Scott reproved this and the subgroup separability of surface groups \cite{Scott78}. More recently, Agol, Long, and Reid \cite{ALR01} proved the Bianchi groups are subgroup separable (see \cite{Agol04B}, \cite{LongReid01}, or \cite{LongReid04B}). In contrast, the mapping class group $\Mod(S)$ of a finite type surface $S$ is known not to be subgroup separable except in a few very special cases (see Section \ref{non}). Nevertheless, it is well-stocked with subgroups of finite index and many interesting subgroups of $\Mod(S)$ are separable.

\begin{theorem} \label{solvable}
Solvable subgroups of $\Mod(S)$ are separable.
\end{theorem}

The proof we give is a variation on the main idea of \cite{RaptisVarsos98}. One of the ingredients readily implies the following.

\begin{theorem}\label{Torelli}
The Torelli group and each term in the Johnson filtration of $\Mod(S)$ is separable.
\end{theorem}

The next theorem gives a potentially large class of subgroups that are separable.
For the statement, set $\Ga_S = \pi_1(S)$, ${\cal X}_n(\Ga_S)$, the $\SL(n;\C)$--character variety, and for a subvariety $V$, $\Stab(V)$ and $\Triv(V)$ the set and pointwise stabilizers of $V$ under the action of $\Mod(S)$. 

\begin{theorem}\label{Variety}
For any number field $k$ and proper $k$--algebraic subvariety $V$ of ${\cal X}_n(\Ga_S)$, the subgroups $\Stab(V)$ and $\Triv(V)$ are separable in $\Mod(S)$.
\end{theorem}

This result is analogous to the linear setting, where for a finitely generated subgroup $\La$ of a $k$--algebraic linear group ${\bf G}$, the subgroup $\La \cap {\bf H}$, for any algebraic subgroup ${\bf H}$ defined over any number field, is separable in $\La$ (see \cite{MargulisSoifer81}, \cite{Long88}, and \cite{Bergeron00}). 

One application of Theorem \ref{Variety} provides separability of the handlebody groups and the Heegaard groups.

\begin{corollary} \label{handle}
The two handlebody groups $\Mod(S,H)$, $\Mod_0(S,H)$ and any Heegaard group $\Mod(S,M^3)$ are separable.
\end{corollary}

Finally, a modest generalization of the proof of the residually finiteness of $\Mod(S)$ \cite{Grossman74} bears our final result.

\begin{theorem} \label{geometric}
The stabilizers of multi-curves, or more generally geometric subgroups of $\Mod(S)$ are separable.
\end{theorem}

\subsection{Acknowledgements}

The authors would like to thank Roger Alperin and Benson Farb for the referral to \cite{RaptisVarsos98}.
We would also like to thank Alan Reid, Jeff Garrett, Bill Goldman, Nick Proudfoot, and Scott Wolpert for useful and interesting conversations.

\section{Definitions and background} \label{backsect}

In this article, we denote a closed orientable surface with genus $g$ and $n$ marked points by $S$,
the surface minus the marked points by $\dot{S}$, and $\Ga_S = \pi_1(\dot{S})$. We define the mapping class group to be the quotient $\Mod(S) = \Diff^+(S)/\Diff_0(S)$, where $\Diff^+(S)$ is the group of orientation preserving diffeomorphisms of $S$ leaving marked points invariant and $\Diff_0(S)$ is the component containing the identity.  It is often useful to note that the natural map from $\Mod(S)$ to $\Out(\Ga_S)$ is injective.  Indeed, when $n = 0$ and $g > 0$, the Dehn-Nielsen Theorem says that this is onto a subgroup of index $2$.

On occasion we will want to consider surfaces-with-boundary and denote these by $\Sigma$ to distinguish them from closed surfaces. In this case we modify the definition of $\Diff^+(\Sigma)$ and demand that boundary components be fixed pointwise.

\subsection{The Johnson filtration}\label{section3}

By a \emph{descending series} on $\Ga_S$, we mean a nested series $\{ C_j \}_{j\geq 0}$ of subnormal subgroups
of $\Ga_S$. When each term $C_j$ is characteristic in $\Ga_S$, we call $C_j$ a \emph{descending characteristic series}. From a descending characteristic series $C_j$ on $\Ga_S$, we obtain a descending series $N_j$ on $\Mod(S)$. Specifically, the canonical set of epimorphisms $\pi_j : \Ga_S \longrightarrow \Ga_S/C_j$ induce a family of homomorphisms $\rho_j :  \Mod(S) \longrightarrow \Out(\Ga_S/C_j)$, and the series $N_j$ is given by $\ker \rho_j$. Following Farb, we refer to the descending series on $\Mod(S)$ induced by the lower central series on $\Ga_S$ as the {\em Johnson filtration}:
\[\Mod(S) \rhd N_1 \rhd N_2 \rhd N_3 \rhd \dots \rhd N_j \rhd \dots.\]
The first non-trivial term $N_1 = {\cal T}$ is usually referred to as the {\em Torelli subgroup} and the second $N_2 = {\cal K}$ the {\em Johnson Kernel}.

\begin{theorem}[Bass-Lubotzky; \cite{BassLubotzky94}]\label{BL}
The Johnson filtration $\{N_j\}_{j> 0}$ on the group $\Mod(S)$ satisfies:
\begin{itemize}
\item[(a)]
\[ \bigcap_j N_j \subset Z(\Mod(S)). \]
\item[(b)]
$N_j/N_{j+1}$ is torsion free for $j>0$.
\end{itemize}
\end{theorem}

\subsection{Solvable subgroups}

Using ideas from the Nielsen-Thurston classification for elements of $\Mod(S)$, Birman, Lubotzky, and McCarthy \cite{BLM83} show that virtually solvable subgroups are virtually abelian (with bounded rank). Ivanov \cite{Ivanov92} strengthened this to the following form.

\begin{theorem}[Birman-Lubotzky-McCarthy, Ivanov]\label{BLM}
There is a finite index subgroup $\Mod'(S) < \Mod(S)$ so that for any virtually solvable subgroup $G < \Mod(S)$, $G \cap \Mod'(S)$ is free abelian with rank at most $3g-3$.
\end{theorem}

\subsection{Multi-curve stabilizers and geometric subgroups}

By a {\em multi-curve}, we mean the isotopy class of a closed embedded $1$--manifold in $\dot{S}$ for which each component is non-peripheral and homotopically essential. $\Mod(S)$ acts on the set of multi-curves, and we denote the stabilizer of a multi-curve $A$ by $\Stab(A)$.

Given a proper subsurface $\Sigma \subset S$, if each component of $\Sigma$ is $\pi_1$--injective, we say that $\Sigma$ is incompressible (which is equivalent to saying that each boundary component is homotopically non-trivial an non-peripheral).
We will often consider subsurfaces as well defined up to isotopy without comment.

An element $\phi \in \Mod(S)$ is called \emph{pure} if for every finite set of $\phi$--invariant connected incompressible subsurfaces $\{ \Sigma_1,...,\Sigma_n \}$ of $S$, each $\Sigma_i$ is invariant.
A subgroup $\La < \Mod(S)$ is called \emph{pure} if every element of $\La$ is pure.
Ivanov \cite{Ivanov92} has shown that there are finite index pure subgroups $\Mod'(S) < \Mod(S)$, and Theorem \ref{BLM} holds for any such subgroup. We note that pure subgroups are, in particular, torsion free.

If $\Sigma \subset S$ is an incompressible subsurface, Paris and Rolfsen \cite{ParisRolfsen00} have proven that the inclusion induces a homomorphism $\Mod(\Sigma) \longrightarrow \Mod(S)$ which is injective, unless two components of the boundary of $\Sigma$ are isotopic in $\dot{S}$. For a general incompressible surface $\Sigma$, the kernel of this homomorphism is contained in the center of $\Mod(\Sigma)$.
Following \cite{ParisRolfsen00}, when $\Sigma$ is incompressible, even if the groups only inject modulo centers, we call the image of $\Mod(\Sigma)$ in $\Mod(S)$ a \emph{geometric subgroup} and denote it by $G(\Sigma)$.

\subsection{Handlebody and Heegaard groups}

In \cite{Masur86}, Masur initiated the study of two subgroups of $\Mod(S)$ called the {\em handlebody groups} defined as follows.
Let $H$ be a handlebody and choose a diffeomorphism $S \longrightarrow \partial H$.
The first handlebody group $\Mod(S,H)$ is the group consisting of those automorphisms of $S$ which extend over $H$.
The second is the subgroup $\Mod_0(S,H)$ of $\Mod(S,H)$ consisting of those elements which induce the identity outer automorphism on $\pi_1(H)$. A similar subgroup studied by Goeritz \cite{Goeritz33}, Powell \cite{Powell80}, and recently Scharlemann \cite{Scharlemann04} is the {\em Heegaard group} $\Mod(S,S^3)$. This is defined by choosing a diffeomorphism from $S$ to a Heegaard surface in $S^3$, then taking the group of automorphisms of $S$ that extend to $S^3$. There is nothing special about $S^3$ in this construction, so for any 3--manifold $M$ and diffeomorphism from $S$ to a Heegaard surface in $M$, we can define the Heegaard group of $(S,M)$, denoted by $\Mod(S,M)$, in an analogous fashion.

We remark that in the present notation, we have suppressed the identification of $S$ with the boundary of the handlebody or Heegaard surface. 

\subsection{Representation and character varieties} \label{repsect}

Given a finitely generated group $\La$ and natural number $n$, the set $\Hom(\La;\SL(n;\C))$ is called the \emph{$\SL(n;\C)$--representation variety}.
This can be equipped with a natural analytic structure (see \cite{Raghunathan72}); in fact, the $\Z$--algebraic structure on $\SL(n;\C)$ provides $\Hom(\La;\SL(n;\C))$ with a $\Z$--algebraic structure.
$\SL(n;\C)$ acts on $\Hom(\La;\SL(n;\C))$ by conjugation and the quotient (in the sense of geometric invariant theory) is a set ${\cal X}_n(\La)$ called the \emph{$\SL(n;\C)$--character variety} which is a $\Z$--defined affine variety (see \cite{MFK94}).

The inner automorphism group action of $\La$ on $\Hom(\La;\SL(n;\C))$ is absorbed by the action of $\SL(n;\C)$. Hence the outer automorphism group $\Out(\La)$ acts naturally on ${\cal X}_n(\La)$ by $\Z$--algebraic automorphisms---see \cite{Magnus80}. For any subvariety $V \subset {\cal X}_n(\La)$ we define
$$\begin{array}{rcl}
\Stab(V) &= &\{ \ga \in \Out(\La) \, : \, \ga(V) \subset V \} \\
\Triv(V) &= &\{ \ga \in \Out(\La) \, : \, \ga|_{V} = \textrm{id}_V \}.
\end{array}$$
If $k$ is a number field, ${\cal O}_k$ its ring of integers, and $V$ is a $k$--defined, then for any finite extension ring $R/{\cal O}_K$, we let $V(R)$ denote the set of $R$--points of $V$.
For an ideal ${\frak p} < R$, reducing the coordinates modulo ${\frak p}$, we obtain a set we denote $V(R/{\frak p})$.  The reduction map (viewed as sets) $r_{{\frak p}} :  V(R) \longrightarrow V(R/{\frak p})$ induces a homomorphism
\[ (r_{{\frak p}})_*  :  \Aut(V(R)) \longrightarrow \Sym(V(R/{\frak p})), \] 
where $\Sym(V(R/{\frak p}))$ is the symmetric group on $V(R/{\frak p})$. Note that since $R/{\frak p}$ is finite, $\Sym(V(R/{\frak p}))$ is a finite group.

\section{Subgroup separability}\label{section2}

\subsection{Separability results}

For convenience we collect the requisite material on separability needed in this article here, referring the reader to the listed references for proofs.

\begin{lemma}[Long-Reid; \cite{LongReid02}] \label{211}
Let $H< K < G$ be groups with $[K:H] <\infty$ and $H$ separable in $G$.
Then $K$ is separable in $G$.
\end{lemma}

\begin{lemma}[Scott; \cite{Scott78}] \label{212}
Let $H,G_0 < G$, with $[G:G_0] < \infty$.
Then $H$ is separable in $G$ if and only if $H \cap G_0$ is separable in $G_0$.
\end{lemma}

\begin{lemma}\label{Surjector}
Let $G$ and $H$ be a pair of groups and $\rho :  G \longrightarrow H$ a surjective homomorphism. If $K<H$ is separable, then $\rho^{-1}(K)$ is separable in $G$.
\end{lemma}

\noindent{\em Proof.}
The lemma is a consequence of the fact that finite index subgroups of $H$ pull back to finite index subgroups of $G$ under $\rho$. \hfill{$\Box$}

\begin{theorem}[Grossman; \cite{Grossman74}]\label{Residual}
$\Mod(S)$ is residually finite.
\end{theorem}

One proof of this theorem uses the fact that $\Ga_S$ satisfies a strong type of residual finiteness: A group $G$ is said to be \emph{conjugacy separable} if for any two non-conjugate elements $x,y \in G$, there is a homomorphism from $G$ to a finite group for which the images of $x$ and $y$ remain non-conjugate. 

\begin{theorem}[Stebe; \cite{Stebe72B}]\label{ConjSep}
$\Ga_S$ is conjugacy separable.
\end{theorem}

\subsection{Lattices and nilpotent groups}

One key result of Segal \cite{Segal90} that we utilize is the following.

\begin{theorem}[Segal]\label{Segal}
The outer automorphism group of a finitely generated nilpotent group is a finite extension of an arithmetic lattice.
\end{theorem}

The work in \cite{McReynolds04A} (see also \cite{AlperinFarb04}) provides the necessary counterpart to this theorem.

\begin{theorem}
Polycyclic subgroups of arithmetic lattices are separable.
\end{theorem}

In combination, these results produce the needed result.

\begin{theorem}\label{Borel}
Polycyclic subgroups of $\Out(\Ga_S/C_j)$ are separable, where $C_j$ is the $j$th term of the lower central series of $\Ga_S$.
\end{theorem}

Using the residual finiteness of $\Out(\Ga_S/C_j)$, it is elementary to deduce Theorem \ref{Torelli}.

\section{Solvable subgroup separability in $\Mod(S)$}\label{section4}

In this section, we prove Theorem \ref{solvable}. Our approach is a variation of the one taken in \cite{RaptisVarsos98} for $\Aut(\Ga_S)$.
By Theorem \ref{BLM}, it suffices to separate abelian subgroups. To this end, we prove the following proposition.

\begin{proposition}\label{Injects}
For every torsion free abelian subgroup $A$ of $\Mod(S)$, there exists $j(A)=j$ such that the quotient
homomorphism $\rho_j :  \Mod(S) \longrightarrow \Mod(S)/N_j$ restricted to $A$ is an injection.
\end{proposition}

\noindent
{\em Proof.}
To begin, set $A_i = \ker\rho_i \cap A = N_i \cap A$. The inclusions $N_{i+1} \leq N_i$ produce inclusions $A_{i+1} \leq A_i$ and so
$\rank_{\Z} A_i \geq \rank_{\Z} A_{i+1}$. As $\{ \rank_{\Z} A_i \}$ is a non-increasing sequence of non-negative integers, there exists $j = j(A) > 0$ such that for all $i,k \geq j$, $\rank_{\Z} A_i = \rank_{\Z} A_k$. By Theorem \ref{BL} (a), $\cap_i A_i = \{ 1 \}$, since $Z(\Mod(S))$ is finite and $A$ is torsion free. Consequently, it suffices to show that the sequence of groups $\{ A_i \}$ is eventually constant. Assuming otherwise, let $i \geq j$ be such that $A_i \neq A_{i+1}$. By our selection of $i$, $\rank_{\Z} A_i = \rank_{\Z} A_{i+1}$, and so $A_i / A_{i+1} < N_i/N_{i+1}$ is a non-trivial finite group. However, this is in opposition with Theorem \ref{BL} (b).\hfill{$\Box$}\\

\noindent
{\em Proof of Theorem \ref{solvable}.}
Fix a pure subgroup $\Mod'(S) < \Mod(S)$ with finite index.
If $A_0$ is any solvable subgroup of $\Mod(S)$, then according to Theorem \ref{BLM} $A = A_0 \cap \Mod'(S)$ is a torsion free abelian group.
By Lemma \ref{212} it suffices to separate $A$ in $\Mod'(S)$.

If we let $N'_i = N_i \cap \Mod'(S)$, then $\cap_i N_i' = \{ 1 \}$; $Z(\Mod(S))$ is finite and so intersects $\Mod'(S)$ trivially.
In an abuse of notation, we write $\rho_i$ for the homomorphism from $\Mod'(S)$ to $\Mod'(S)/N_i'$, and note that Proposition \ref{Injects} still holds for some $j = j(A)$.

For any $i \geq j$, $\rho_i^{-1}(\rho_i(A)) = A N'_i = A \ltimes N'_i$ since $\rho_i$ restricted to $A$ is injective.
An application of Theorem \ref{Borel} implies $\rho_i(A)$ is separable. Hence, by Lemma \ref{Surjector}, $\rho_i^{-1}(\rho_i(A))$ is separable. The inclusion of $A \ltimes N'_k \hookrightarrow A \ltimes N'_i$ for $k > i$ respects the semidirect product structure. This in combination with $\cap_i N'_i = \{ 1 \}$ yields
\[ \bigcap_{k \geq i} \rho_k^{-1}(\rho_k(A)) = \bigcap_{k \geq i} A \ltimes N^\prime_k = A \ltimes \bigcap_{k \geq i} N_k^\prime = A. \] 
Since $\rho_k^{-1}(\rho_k(A))$ is separable, it is expressible as an intersection of finite index subgroups, hence so is $A$.\hfill{$\Box$}\\

This argument along with the theorem for $\Out(F_n)$ analogous to Theorem \ref{BL}---see \cite{BassLubotzky94}---yields

\begin{theorem}\label{Free}
Solvable subgroups of $\Out(F_n)$ are separable.
\end{theorem}

The proof of Theorem \ref{Free} also requires an analog of Proposition \ref{Injects} for polycyclic subgroups of $\Out(F_n)$ which is achieved by considering Hirsch length as opposed to rank. That this yields all solvable subgroups of $\Out(F_n)$ follows from the fact that solvable subgroups are virtually polycyclic---see \cite{BassLubotzky94}. 

\section{Subvariety stabilizer separability}

To prove Theorem \ref{Variety} we would like to use the reduction maps from Section \ref{repsect} to construct finite index subgroups of $\Mod(S)$. This, in turn, requires the existence of sufficiently many algebraic points on the variety $V$. The starting point for our proof of Theorem \ref{Variety} is thus the following consequence of Hilbert's Nullstellensatz.

\begin{proposition}\label{Density}
Let $V$ be a $k$--algebraic variety for a number field $k$. Then
$$V(\overline{k}) \quad = \bigcup_{\begin{array}{c} k<K \\ | K:k | < \infty \end{array}}V(K)$$
is Zariski dense in $V$, where $\overline{k}$ is the algebraic closure of $k$.
\end{proposition}

Given a finite extension $K/k$ and $q \in V(K)$, it follows that $q \in V(R)$, for some finite extension ring $R/{\cal O}_K$. For example, if $m$ is the product of the denominators occurring in the coordinates of $q$, one can take $R={\cal O}_K[1/m]$. Coupled with Proposition \ref{Density}, this produces a plethora of points algebraically defined in the variety $V$.\\

\noindent{\em Proof of Theorem \ref{Variety}.}
Let $V \subset {\cal X}_n(\Ga_S)$ be a $k$--algebraic subvariety, and let ${\frak a} \subset {\frak a}_V \subset \C[\bf{T}]$ be the associated ideals for ${\cal X}_n(\Ga_S)$ and $V$, respectively.
For $\ga \in \Mod(S) \setminus \Stab(V)$, there exists $q_0 \in V$ such that $\ga(q_0) \notin V$.
Since $V(\overline{k})$ is Zariski dense in $V$, we can assume $q_0 \in V(K)$, for some finite extension $K/k$. Indeed, from the discussion above, $q_0 \in V(R)$ for some finite extension ring $R/{\cal O}_K$. Next, select a generating set $f_1,\dots,f_r \in R[\bf{T}]$ for ${\frak a}_V$. It follows that $f_j(q_0),f_j(\ga(q_0)) \in R$ for all $j$. By assumption, $\ga(q_0) \notin V$, and so $f_j(\ga(q_0)) \ne 0$ for some generator $f_j$.
Thus, we can select an ideal ${\frak p}$ of $R$ such that $f_j(\ga(q_0)) \ne 0 \mod {\frak p}$.
Since the points of $V(R/{\frak p})$ are precisely those of the form $r_{\frak p}(q)$, with $f_i(q) = 0 \mod {\frak p}$ for each $1 \leq i \leq r$, it follows that $r_{\frak p}(\ga(q_0)) \notin V(R/{\frak p})$.
Since $(r_{\frak p})_*(\ga)(r_{\frak p}(q_0)) =r_{\frak p}(\ga(q_0))$, we have $(r_{\frak p})_*(\ga) \notin \Stab(V/{\frak p})$. On the other hand, because $r_{\frak p}(V(R)) = V(R/{\frak p})$, we see that $(r_{\frak p})_*(\Stab(V))$ is contained in $\Stab(V(R/{\frak p}))$. Therefore, $(r_{\frak p})_*^{-1}(\Stab(V(R/{\frak p})))$ separates $\ga$ from $\Stab(V)$.

Next let $\ga \in \Mod(S) \setminus \Triv(V)$.
As above, we may select $q_0 \in V(R)$, for some finite extension $R/{\cal O}_K$, such that $\ga(q_0) \ne q_0$.
The separation of $\Triv(V)$ from $\ga$ is obtained in a similar fashion to the above by selecting an ideal ${\frak p}$ of $R$ such that $\ga(q_0) \ne q_0 \mod {\frak p}$. Specifically, with such an ideal, $(r_{\frak p})_*^{-1}(\Triv(V(R/{\frak p})))$ separates $\ga$ and $\Triv(V)$. \hfill{$\Box$}\\

As mentioned in the introduction, separability of the handlebody groups and Heegaard groups follows from Theorem \ref{Variety}.\\

\noindent
{\em Proof of Corollary \ref{handle}.}
A Heegaard group is the intersection of two conjugates of the handlebody group $\Mod(S,H)$, so it suffice to prove the corollary for $\Mod(S,H)$ and $\Mod_0(S,H)$.

We note that the embedding $i :  S \longrightarrow H$ induces an inclusion
$$i^* : {\cal X}_2(\pi_1(H)) \longrightarrow {\cal X}_2(\Ga_S).$$
The image under $i^*$, denoted by $V_H$, is $\Z$--defined. Specifically, its obtained by declaring certain elements in $\Ga_S$---those that bound disks in $H$---to be trivial.
Clearly, $\Mod(S,H) < \Stab(V_H)$.
If $\phi \in \Mod(S) \setminus \Mod(S,H)$ then it takes a simple closed curve that bounds a disk in $H$ to one which does not bound a disk, so by Dehn's Lemma \cite{Rolfsen76} it is not homotopically trivial in $H$.
Since there exist faithful representations of $\pi_1(H)$ to $\SL(2;\C)$, it follows that $\phi \notin \Stab(V_H)$, hence $\Stab(V_H) = \Mod(S,H)$, and $\Mod(S,H)$ is separable by Theorem \ref{Variety}.
Similarly, $\Mod_0(S,H) = \Triv(V_H)$ and so is separable by Theorem \ref{Variety}.\hfill{$\Box$}

\begin{remark}
The idea of using algebraic actions of groups on varieties to deduce residual properties is not new. Bass and Lubotzky \cite{BassLubotzky83} used this to produce another proof of the residual finiteness of $\Mod(S)$.
\end{remark}

\section{Geometric subgroups}\label{section5}

The next theorem easily implies Theorem \ref{geometric}.

\begin{theorem}\label{fstab}
If $\Delta$ is any finite set of conjugacy classes in $\Ga_S$, then the group $\Stab(\Delta)$ is separable.
\end{theorem}

\noindent
{\em Proof.}
Write the elements of $\Delta$ as $\Delta = \{a_1, \dots, a_n\}$ .
Let $\phi \in \Mod(S) \setminus \Stab(\Delta)$ and $a_i$ be such that $\phi(a_i) = x$ and $x \neq a_j$ for any $j$.
By Theorem \ref{ConjSep}, there exists a finite group $F$ and epimorphism
$\psi : \Ga_S \longrightarrow F$ for which $\psi(x)$ is not conjugate to $\psi(a_i)$ for any $i$.
Without loss of generality, we may assume that $\ker(\psi)$ is characteristic (if not, replace $F$ by the quotient of $\Ga_S$ by the characteristic core of $\ker(\psi)$). The homomorphism $\psi$ induces a homomorphism
\[ \psi_* : \Out(\Ga_S) \longrightarrow \Out(F). \] 
The group $\Out(F)$ acts on the set of conjugacy classes of $F$ and we can consider the subgroup $\Stab(\psi(\Delta))$. Visibly, $\Stab(\Delta) \subset \psi^{-1}(\Stab(\psi(\Delta)))$ and does not contain $\psi_*(\phi)$ since $\psi(x)$ is not in any of the conjugacy classes of $\psi(\Delta)$. \hfill{$\Box$}\\

\noindent
{\em Proof of Theorem \ref{geometric}.}
The proof for multi-curve stabilizers is immediate from Theorem \ref{fstab}. Suppose $\Sigma \subset S$ is an incompressible subsurface. To separate $G(\Sigma)$ we show that for any finite index pure subgroup $\Mod'(S) < \Mod(S)$
\[ G(\Sigma) \cap \Mod'(S) = \Stab(\Delta) \cap \Mod'(S) \]
for some multi-curve $\Delta$.
Then, by Lemma \ref{212} and Theorem \ref{fstab}, $G(\Sigma)$ will be separable.

We construct $\Delta$ as follows. Let $\partial \Sigma$ denote the union of the conjugacy classes in $\Ga_S$ representing the boundary components of $\Sigma$ and $R_1,...,R_k$ denote those components of the complementary subsurface $\overline{S \setminus \Sigma}$ that are not homeomorphic to annuli or pairs of pants. For each $R_i$, let $\alpha_i,\beta_i$ be a pair of simple closed curves which bind $R_i$ (meaning every other non-peripheral essential closed curve on $R_i$ intersects one of $\alpha_i$ or $\beta_i$). Any automorphism of $R_i$ that fixes both $\alpha_i$ and $\beta_i$ must have finite order (up to Dehn twisting in curves parallel to the boundary)---compare \cite{Kerckhoff83}. Finally, set $\Delta = \partial \Sigma \cup \{\alpha_1,\beta_1,...,\alpha_k,\beta_k\}$.

If $\phi$ is any pure automorphism which preserves $\Delta$, it must be the identity on $\overline{S \setminus \Sigma}$.
This is true for each $R_i$ by the discussion in the previous paragraph, the definition of pure automorphism, and the fact that Dehn twisting in the boundary of $\overline{S \setminus \Sigma}$ can be absorbed into $\Sigma$.
On the other hand, since the automorphism group of an annulus or pair of pants is finite, modulo Dehn twisting in the boundary, the purity assumption shows that it holds for these components as well.
Thus $\Stab(\Delta) \cap \Mod'(S) < G(\Sigma) \cap \Mod'(S)$.
Since being in $\Stab(\Delta)$ puts absolutely no constraint on the restriction to $\Sigma$, the other inclusion evidently holds.\hfill{$\Box$}

\begin{remark}
It is worth mentioning that conjugacy separability is equivalent to the quotient topology on the set of conjugacy classes induced by the profinite topology on $\Ga_S$ being Hausdorff. Consequently, finite sets are closed. Using the proof of Theorem \ref{fstab} as a guide, one can prove:
\end{remark}

\begin{theorem}
The stabilizer in $\Mod(S)$ of any closed set of conjugacy classes in $\Ga_S$ is separable.
\end{theorem}

\section{Appendix: $\Mod(S)$ is not subgroup separable}\label{non}

For completeness, in this section we prove that $\Mod(S)$ is not subgroup separable. The idea
is to show that $\Mod(S)$ contains a subgroup which is not subgroup separable. To begin, whenever $S$ contains two disjoint incompressible subsurfaces more complicated than a pair of pants, $\Mod(S)$ contains an isomorphic copy of $F_2 \times F_2$, and so cannot be subgroup separable. In the remaining cases, $\Mod(S)$ is virtually free, and hence obviously separable, or else $S$ is a torus with two marked points or a sphere with 5 marked points. In these cases, $\Mod(S)$ contains the fundamental group of every 3--manifold fibered over the circle with fiber a once punctured torus or 4--punctured sphere, respectively. It is well known that there are such 3--manifold groups (with reducible monodromy) which are not subgroup separable (see \cite{BKS87}). Therefore, the mapping class groups in these special cases are not subgroup separable.



\def\lfhook#1{\setbox0=\hbox{#1}{\ooalign{\hidewidth
  \lower1.5ex\hbox{'}\hidewidth\crcr\unhbox0}}} \def\cprime{$'$}
  \def\cprime{$'$}
\providecommand{\bysame}{\leavevmode\hbox to3em{\hrulefill}\thinspace}
\providecommand{\MR}{\relax\ifhmode\unskip\space\fi MR }
\providecommand{\MRhref}[2]{%
  \href{http://www.ams.org/mathscinet-getitem?mr=#1}{#2}
}
\providecommand{\href}[2]{#2}


\noindent
Department of Mathematics, Columbia University, 2990 Broadway MC 4448, New York, NY 10027-6902\\
email: {\tt clein@math.columbia.edu}\\

\noindent
Department of Mathematics, The University of Texas at Austin, 1 University Station C1200, Austin, TX 78712-0257\\
email: {\tt dmcreyn@math.utexas.edu}

\end{document}